\documentclass[10pt]{article}

\usepackage[english]{babel}
\usepackage{ae}
\usepackage{graphicx}
\usepackage{tikz}
\usetikzlibrary{matrix}
\usepackage{amsmath, amsfonts, amssymb}
\usepackage{url}
\usepackage{authblk}

\providecommand{\abs}[1]{\lvert#1\rvert}

\providecommand{\norm}[1]{\lVert#1\rVert}

\newcommand{\num}[1]{\mathbb{#1}}
\newcommand{\Reals}{\num{R}}
\newcommand{\cali}[1]{\mathcal{#1}}

\def\vec#1{\mathchoice{\mbox{\boldmath$\displaystyle#1$}}
{\mbox{\boldmath$\textstyle#1$}}
{\mbox{\boldmath$\scriptstyle#1$}}
{\mbox{\boldmath$\scriptscriptstyle#1$}}}

\begin{document}

\title{Optical Flow on Evolving Surfaces with an Application to the Analysis of\\ 4D Microscopy Data}

\author[1]{Clemens Kirisits}
\author[1]{Lukas F. Lang}
\author[1,2]{Otmar Scherzer}
\affil[1]{\footnotesize Computational Science Center, University of Vienna, Nordbergstr.\ 15, 1090 Vienna, Austria}
\affil[2]{Radon Institute of Computational and Applied Mathematics, Austrian Academy of Sciences, Altenberger Str.\ 69, 4040 Linz, Austria}

\maketitle

\begin{abstract}
\noindent
We extend the concept of optical flow to a dynamic non-Euclidean setting. Optical flow is traditionally computed from a sequence of flat images. It is the purpose of this paper to introduce variational motion estimation for images that are defined on an evolving surface. Volumetric microscopy images depicting a live zebrafish embryo serve as both biological motivation and test data.
\end{abstract}

\noindent
\textbf{Keywords: }Computer Vision, biomedical imaging, optical flow, variational methods, evolving surfaces, zebrafish, laser-scanning microscopy.

\section{Introduction}
Advances in laser-scanning microscopy and fluorescent protein technology have increased resolution of microscopy imaging up to a single cell level~\cite{MegFra03}. They allow for four-dimensional (volumetric time-lapse) imaging of living organisms and shed light on cellular processes during early embryonic development. Understanding cellular development often requires estimation and analysis of cell motion. However, the amount of data captured is tremendous and therefore manual analysis is not an option.

\begin{figure}[t]
	\centering
\includegraphics[width=0.45\textwidth]{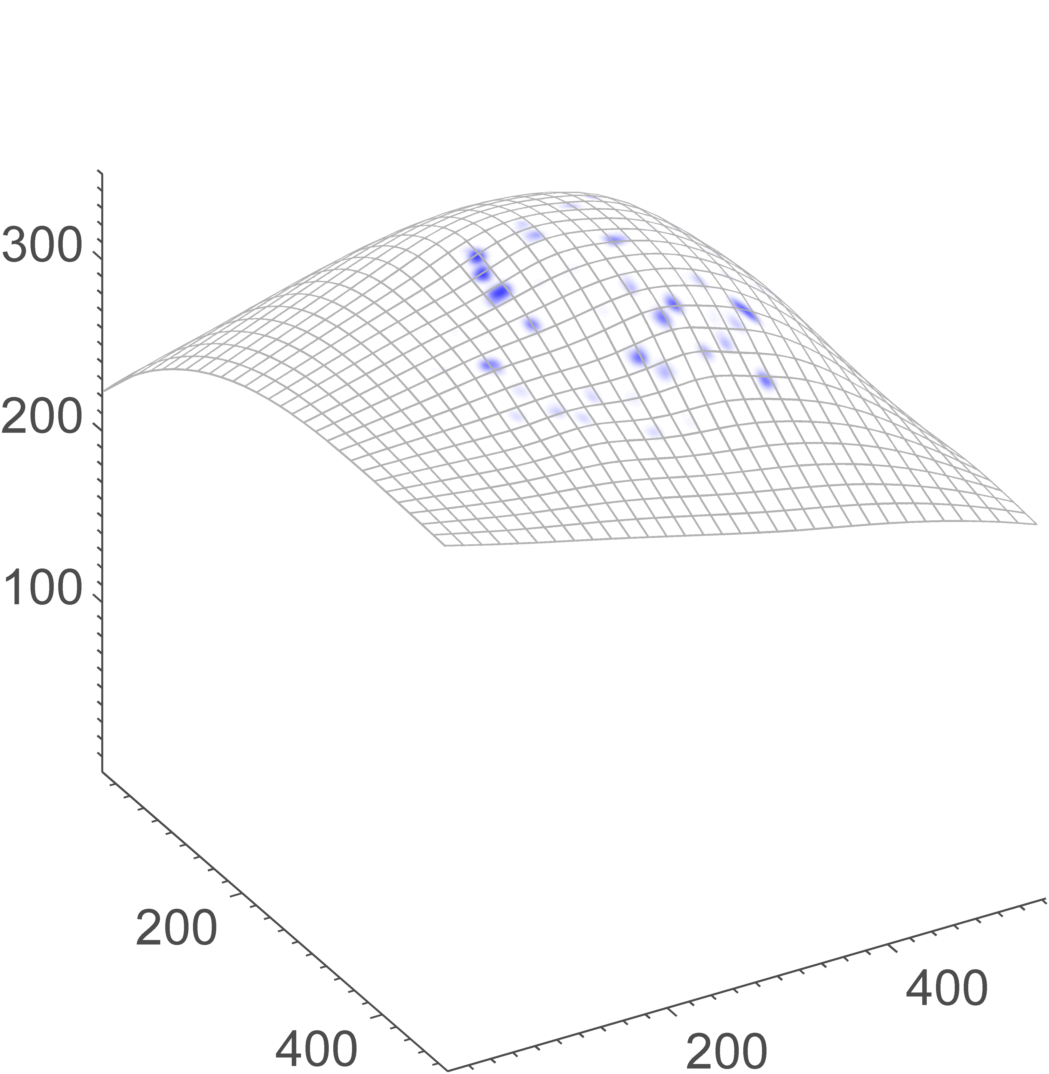}
	\hfill
\includegraphics[width=0.45\textwidth]{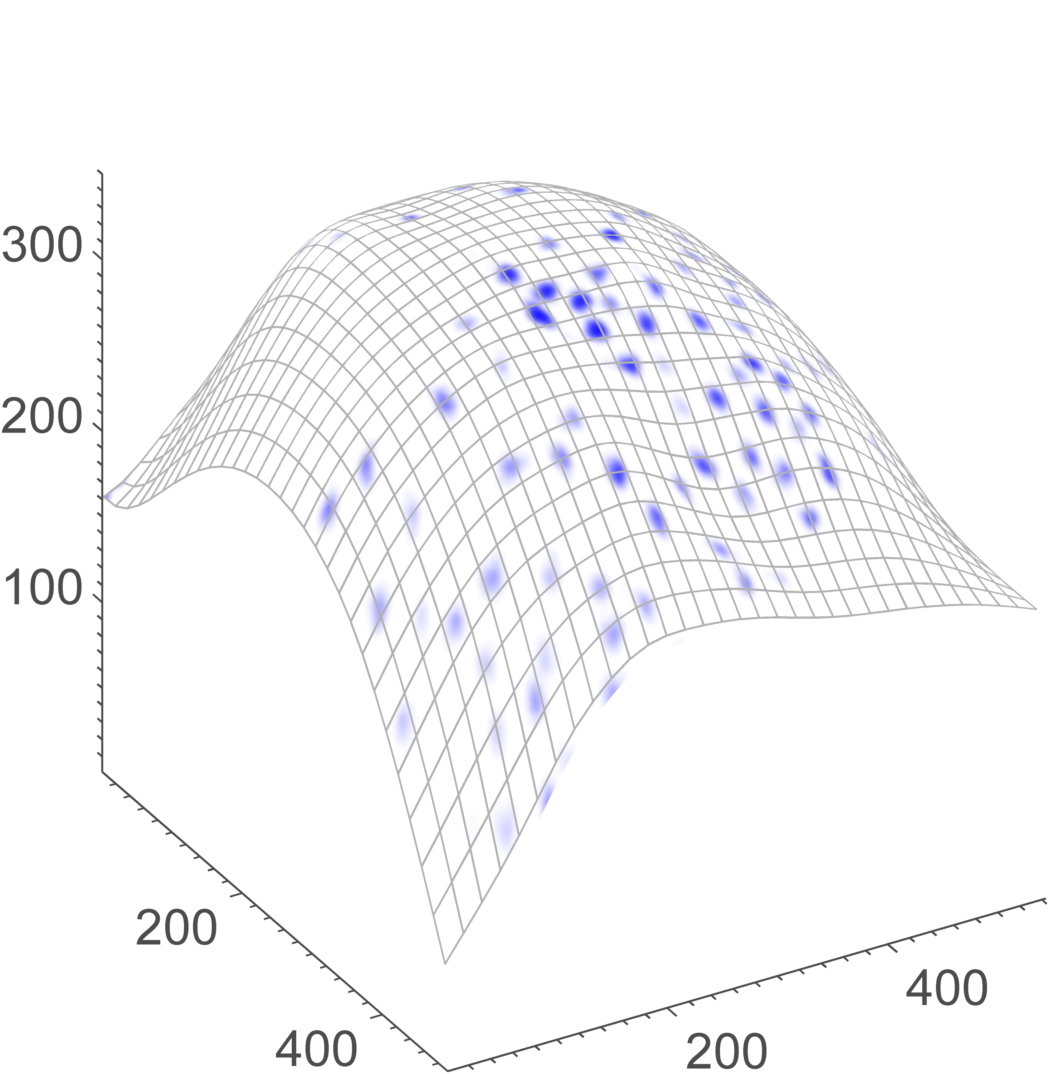} \\
	\medskip
\includegraphics[width=0.45\textwidth]{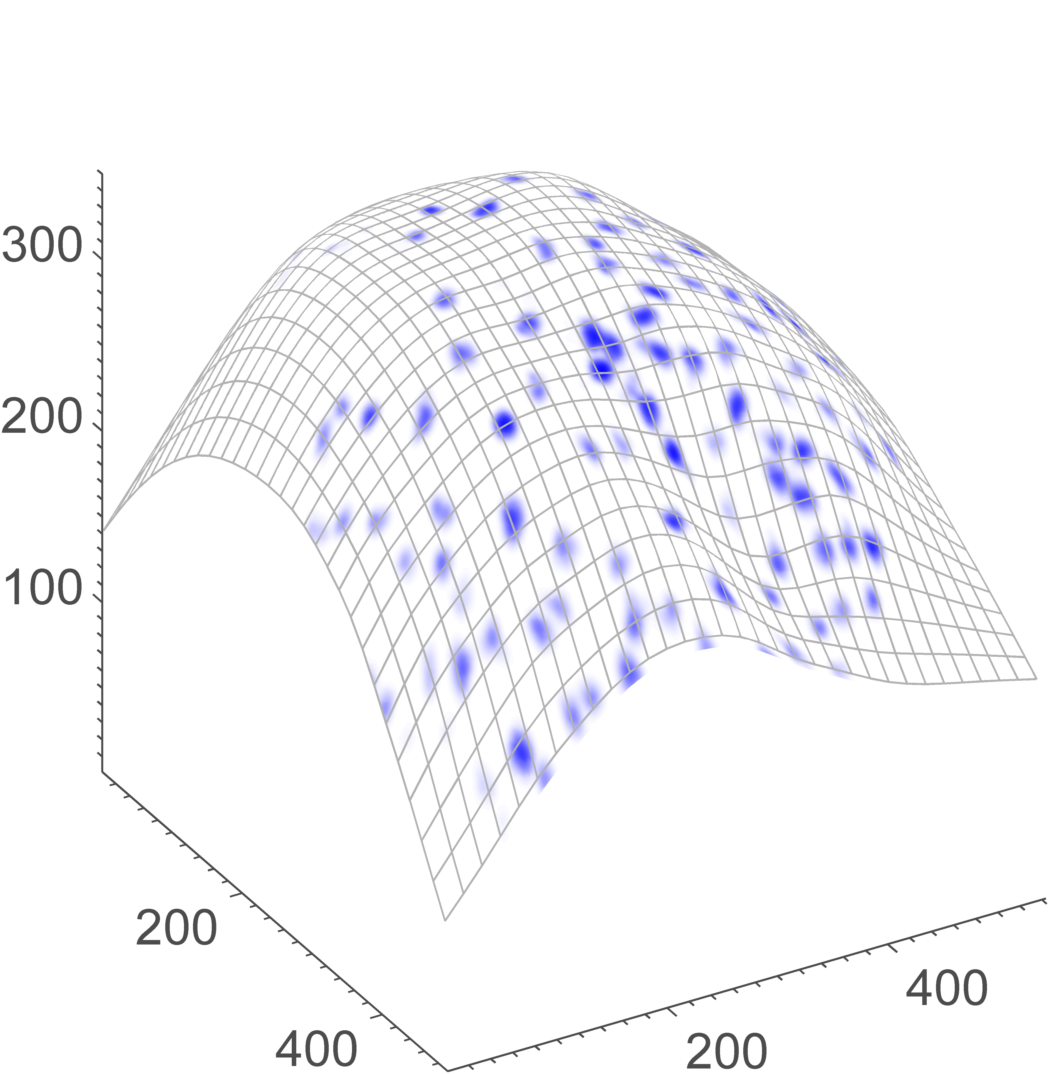}
	\hfill	
\includegraphics[width=0.45\textwidth]{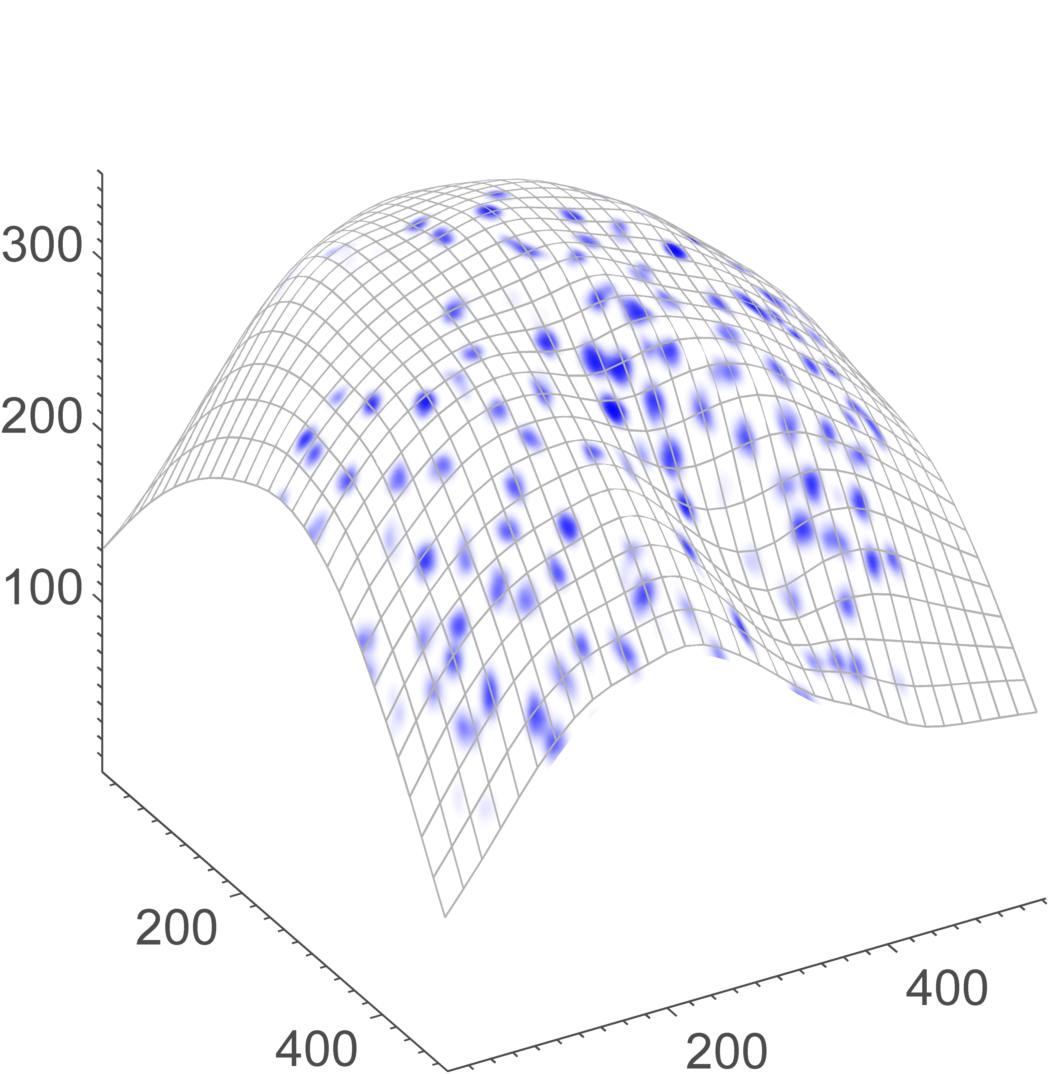}
	\caption{Sequence of embryonic zebrafish images. The curved mesh represents a section of the yolk's surface. Depicted are frames no.\ 30, 45, 55, and 60 of the entire sequence. All dimensions are in micrometer ($\mu$m). See Sec.~\ref{sec:data} for more details on the microscopy data.}
	\label{fig:embryo}
\end{figure}

The specific biological motivation for this work is to understand the motion and division behaviour of fluorescently labelled endodermal cells of a zebrafish embryo. The marked cells develop on the surface of the embryo's yolk, where they form a non-contiguous monolayer \cite{WarNus99}. Loosely speaking, they only sit next to each other but not on top of each other. Moreover, the yolk deforms over time; see Fig.\ \ref{fig:embryo}.

We take these biological facts into account and restrict our attention to the analysis of cell motion on the yolk's surface. With this approach it is possible to reduce the amount of data by one space dimension. The resulting problem consists in the estimation of motion of brightness patterns that are restricted to an itself moving surface. We approach this problem by adapting the classical concept of optical flow to the present setting, where the image domain is both non-Euclidean and dynamic. Note that due to the monolayer structure cell occlusions cannot occur. This makes the optical flow field a more reliable approximation to the true motion field.

Our contributions in the field of optical flow are as follows. First, we formulate the optical flow problem on an evolving two-dimensional manifold and give two equivalent ways of linearising the brightness constancy assumption (Secs.~\ref{sec:bca} and \ref{sec:ofc}). One uses a parametrisation of the evolving surface, the other one is parameter-independent. Second, we use a generalisation of the Horn-Schunck model to regularise the optical flow field (Sec.~\ref{sec:regularisation}). For a given global parametrisation of the evolving surface, we solve the associated Euler-Lagrange equations in the parameter domain with a finite difference scheme (Sec.~\ref{sec:num}). Finally, we apply this technique to obtain qualitative results from the afore-mentioned zebrafish data (Sec.~\ref{sec:exp}). Our experiments show that the optical flow is an appropriate tool for analysing these data. It is capable of estimating global trends as well as individual cell movements and, in particular, it is able to indicate cell division events.

\subsubsection{Related work.}
Optical flow is the apparent motion in a sequence of images. Its estimation is a key problem in Computer Vision. Horn and Schunck~\cite{HorSchu81} were the first to propose a variational approach assuming constant brightness of moving points and spatial smoothness of the velocity field. Since then, a vast number of modifications has been developed. See \cite{BakSchaLewRotBla11} for a recent survey.

Miura \cite{Miu05} observed that until 2005 optical flow has been mostly disregarded as a method for motion extraction in cell biological data. Since then, a few articles have explored this direction: Melani et al.~\cite{MelCamLomRizVer07} and Huben\'{y} et al.~\cite{HubUlmMat07} extended variational optical flow methods to volumetric images to obtain 3D displacement fields. In the former article, the resulting algorithm is also applied to zebrafish microscopy data. Quelhas et al.~\cite{QueMenCam10} use optical flow to detect cell divisions in a live plant root. However, they work with 2D (plus time) data only. Therefore, their approach suffers from errors caused by 3D off-plane motion.

Clearly, certain natural scenarios are more accurately described by a velocity field on a non-flat surface rather than on a flat domain. With applications to robot vision, Imiya~et~al.~\cite{ImiSugTorMoc05,TorImiSugMoc05} considered optical flow for spherical images. In a more general setting, Lef\`{e}vre and Baillet~\cite{LefBai08} extended the Horn-Schunck method to 2-Riemannian manifolds and showed well-posedness. They solve the numerical problem with finite elements on a surface triangulation. In all of the above works the underlying imaging surface is fixed over time, while in this paper it is not.
\section{Optical Flow on Evolving Surfaces} \label{sec:model}
\subsection{Brightness Constancy}\label{sec:bca}
Let $\cali{M}_t \subset \Reals^3$, $t \in I = [0,T)$, be a compact smooth two-dimensional manifold evolving smoothly over time. We assume the velocity to be unknown. Moreover, denote by $\tilde{f}$ a scalar time-dependent quantity defined on the surface
\begin{equation*}
	\tilde{f} \colon \bigcup_{t \in I} \left( \cali{M}_t \times \{t\} \right) \to \Reals.
\end{equation*}
We begin with a Lagrangian specification of the optical flow field. That is, for every starting point $\vec{x}_0 \in \cali{M}_0$ we seek a trajectory where the data $\tilde{f}$ are conserved. More precisely, we want to find a function
\begin{equation*}
	\gamma \colon \cali{M}_0 \times I \to \bigcup_{t \in I} \cali{M}_t, 
\end{equation*}
such that
\begin{enumerate}
	\item	$\gamma(\vec{x}_0,t) \in \cali{M}_t$ for all $t \in I$, for all $\vec{x}_0 \in \mathcal{M}_0$,
	\item	$\gamma(\cdot,t)$ is a diffeomorphism between $\cali{M}_0$ and $\cali{M}_t$ for all $t \in I$,
	\item	$\gamma(\cdot,0) = \mathrm{Id}_{\cali{M}_0}$,
\end{enumerate}
is fulfilled and which satisfies a ``brightness" constancy assumption (BCA)
\begin{equation}\label{eq:bca}
	\tilde{f}(\vec{x}_0,0) = \tilde{f}(\gamma(\vec{x}_0,t),t), \enspace \mbox{for all } (\vec{x}_0,t) \in \cali{M}_0 \times I.
\end{equation}

In classical optical flow computations it is common practice to linearise the BCA by taking its time derivative and to solve the resulting equation for the Eulerian unknown $\dot{\gamma}$.\footnote{To simplify expressions we use Newton's notation for those time derivatives that correspond to actual velocities, for example $\dot\gamma = \partial_t \gamma$.} We also take this route, but differentiation of $\tilde{f}$ is more involved. Observe, for example, that for an arbitrary $t_0 \in I $ and $\vec{x} \in \cali{M}_{t_0}$ the usual partial derivative
\begin{equation*}
	\partial_t \tilde{f} (\vec{x},t_0) = \lim_{h \to 0} \frac{1}{h} \left( \tilde{f}(\vec{x},t_0+h) - \tilde{f}(\vec{x},t_0) \right)
\end{equation*}
is not well-defined, simply because, in general, $\vec{x}$ is not an element of $\cali{M}_{t_0+h}$ for all $h \neq 0$.

In the next section we linearise \eqref{eq:bca} in two different ways. First, we use a global parametrisation to pull the data back to a fixed reference domain and linearise afterwards. In our second approach we borrow some notions from continuum mechanics \cite{CerFriGur05} to directly linearise \eqref{eq:bca}.

\subsection{Linearisation}\label{sec:ofc}
\paragraph{Linearisation after pull-back.}
Let $\Omega \in \Reals^2$ be a compact domain and
\begin{equation*}
	\vec{x} \colon \Omega \times I \to \Reals^3, \enspace (x_1, x_2, t) = (x,t)  \mapsto \vec{x}(x,t) \in \cali{M}_t
\end{equation*}
be a parametrisation of the evolving surface. Denote by $f$ the coordinate representation of $\tilde{f}$, that is,
\begin{equation}\label{eq:data}
	f(x,t) = \tilde{f}(\vec{x}(x,t),t)
\end{equation}
and let
\begin{equation*}
	\beta \colon \Omega \times I \to \Omega
\end{equation*}
be the coordinate counterpart of $\gamma$. This means, if we let $\vec{x}_0 = \vec{x}(x_0,0)$, then $\beta(x_0,t)$ gives the coordinates of $\gamma(\vec{x}_0,t) \in \cali{M}_t$ in $\Omega$ (see Fig.~\ref{fig:cd}). In other words, we have the identity
\begin{equation}\label{eq:betagamma}
	\gamma(\vec{x}(x_0,0),t) = \vec{x}(\beta(x_0,t),t), \enspace \mbox{ for all } (x_0,t) \in \Omega \times I.
\end{equation}
Now, from \eqref{eq:bca}, \eqref{eq:data} and \eqref{eq:betagamma} we get
\begin{align*}
	f(x_0,0)	&=	\tilde{f}(\vec{x}_0,0) \\
				&=	\tilde{f}(\gamma(\vec{x}_0,t),t) \\
				&=	\tilde{f}(\vec{x}(\beta(x_0,t),t),t) \\
				&=	f(\beta(x_0,t),t),
\end{align*}
which is a coordinate version of the BCA. After differentiation with respect to $t$ it becomes
\begin{equation}\label{eq:ofc}
	\nabla^2 f \cdot \dot{\beta} + \partial_t f = 0,
\end{equation}
where $\nabla^2 = (\partial_1, \partial_2)^\top$ is the two-dimensional spatial gradient. Note that the last equation is nothing but the classical optical flow constraint (OFC) for Euclidean data $f$ and a displacement field $\dot{\beta}$.

\begin{figure}[t]
	\centering
	\begin{tikzpicture}
		\matrix (m) [matrix of math nodes,row sep=3em,column sep=4em,minimum width=2em]{
	    	\Omega & \Omega \\
	    	\mathcal{M}_0 & \mathcal{M}_t \\};
		\path[-stealth]
	    	(m-1-1) edge node [left]	{$\vec{x}(\cdot,0)$}	(m-2-1)
	    			edge node [above]	{$\beta(\cdot,t)$}		(m-1-2)
	    	(m-2-1) edge node [above]	{$\gamma(\cdot,t)$}		(m-2-2)
	    	(m-1-2) edge node [right]	{$\vec{x}(\cdot,t)$}	(m-2-2);
	\end{tikzpicture}
	\caption{Commutative diagram describing the relation between unknowns $\beta$ and $\gamma$.}
	\label{fig:cd}
\end{figure}

\paragraph{Direct linearisation.}
We turn to our second derivation. While, as pointed out above, the partial derivative $\partial_t\tilde{f}$ is undefined in general, it does make sense to differentiate $\tilde{f}$ following the surface movement. Let $\vec{y}$ be a point on $\cali{M}_{t_0}$ and $\xi \colon t \mapsto \xi(t) \in \cali{M}_t$ an arbitrary smooth trajectory through the evolving surface satisfying $\xi(t_0) = \vec{y}$. Now we can compute
\begin{equation*}
	\left. \frac{\mathrm{d}}{\mathrm{d}t} \tilde{f}(\xi(t),t)\right|_{t=t_0} = \lim_{h \to 0} \frac{1}{h} \left( \tilde{f}(\xi(t_0+h),t_0+h) - \tilde{f}(\vec{y},t_0) \right)
\end{equation*}
to obtain a valid derivative of $\tilde{f}$. Since this time derivative only depends on the vector $\vec{v} = \dot{\xi}(t_0)$, we denote it by $\mathrm{d}^{\vec{v}}_t \tilde{f}$. A natural candidate for a trajectory along which to differentiate is given by the parametrisation $\xi(t) = \vec{x}(x,t)$. Another possible choice would be a trajectory that is normal to $\cali{M}_{t_0}$. The resulting normal time derivative is accordingly denoted by $\mathrm{d}^{\vec{n}}_t\tilde{f}$.

Finally, we also need the surface gradient $\nabla_{\cali{M}} \tilde{f}$. If $F$ is a smooth extension of $\tilde{f}$ to an open neighbourhood of $\vec{y} \in \cali{M}_{t_0}$ in $\Reals^3$, then the surface gradient of $F$ at $\vec{y}$ is defined as the projection of the three-dimensional spatial gradient $\nabla^3 F$ onto the tangent plane to $\cali{M}_{t_0}$
\begin{equation*}
	\nabla_{\cali{M}} F = \nabla^3 F - (\nabla^3 F \cdot \hat{\vec{n}})\hat{\vec{n}},
\end{equation*}
where $\hat{\vec{n}}$ is the unit normal to $\cali{M}_{t_0}$. The surface gradient only depends on the values of $F$ on the surface; see e.g.\ \cite[p.~389]{GilTru01}. Thus, $\nabla_{\cali{M}} \tilde{f} = \nabla_{\cali{M}} F$ is well-defined.

The spatial and temporal derivatives of $\tilde{f}$ introduced above are related in a simple way. As shown in \cite{CerFriGur05}, they satisfy the equality
\begin{equation}\label{eq:gurtin}
	\begin{aligned}
	\mathrm{d}^{\dot{\vec{x}}}_t \tilde{f}
		&= \nabla_{\cali{M}} \tilde{f} \cdot \dot{\vec{x}} + \mathrm{d}^{\vec{n}}_t\tilde{f} \\
		&= \nabla_{\cali{M}} \tilde{f} \cdot \dot{\vec{x}}_{\mathrm{tan}} + \mathrm{d}^{\vec{n}}_t\tilde{f},
	\end{aligned}
\end{equation}
where $\dot{\vec{x}}_{\mathrm{tan}}$ is the tangential surface velocity, that is, the projection of $\dot{\vec{x}}$ onto the tangent plane to $\cali{M}_{t_0}$. This decomposition of $\mathrm{d}^{\dot{\vec{x}}}_t \tilde{f}$ into normal and tangential components is clearly valid for any trajectory in place of $\vec{x}$, and therefore in particular for the unknown $\gamma$. This means we can use \eqref{eq:gurtin} in order to differentiate the BCA \eqref{eq:bca} with respect to $t$. The resulting OFC reads
\begin{equation}\label{eq:ofc2}
	\nabla_{\cali{M}} \tilde{f} \cdot \dot{\gamma}_\mathrm{tan} + \mathrm{d}^{\vec{n}}_t\tilde{f} = 0.
\end{equation}

\paragraph{Discussion.}
We conclude this section with a brief comparison of the two OFCs derived above. We start by showing how to obtain \eqref{eq:ofc} from \eqref{eq:ofc2} and vice versa. To this end we again assume the existence of a global parametrisation and rewrite all quantities in \eqref{eq:ofc2} in terms of $\vec{x}$. First observe that, by \eqref{eq:betagamma}, the velocity of $\gamma$ equals the surface velocity $\dot{\vec{x}}$ plus a purely tangential component
\begin{equation*}\label{eq:velocity}
	\dot{\gamma} = \dot{\vec{x}} + J \dot{\beta},
\end{equation*}
where $J = (\partial_1 \vec{x} \; \partial_2 \vec{x})$ is the Jacobian matrix of $\vec{x}$ with respect to $x$. On the other hand, by \eqref{eq:gurtin}, the normal time derivative is equal to the time derivative of $\tilde{f}$ following $\vec{x}$ minus its tangential component
\begin{equation*}
	\mathrm{d}^{\vec{n}}_t\tilde{f} = \mathrm{d}^{\dot{\vec{x}}}_t \tilde{f} - \nabla_{\cali{M}} \tilde{f} \cdot \dot{\vec{x}}.
\end{equation*}
Using the last two equations to rewrite the left-hand side of \eqref{eq:ofc2} yields
\begin{align*}
	\nabla_{\cali{M}} \tilde{f} \cdot \dot{\gamma} + \mathrm{d}^{\vec{n}}_t\tilde{f}
		&=	\nabla_{\cali{M}} \tilde{f} \cdot \left(\dot{\vec{x}} + J \dot{\beta}\right) + \mathrm{d}^{\dot{\vec{x}}}_t \tilde{f} - \nabla_{\cali{M}} \tilde{f} \cdot \dot{\vec{x}} \\
		&=	\nabla_{\cali{M}} \tilde{f} \cdot J \dot{\beta} + \mathrm{d}^{\dot{\vec{x}}}_t \tilde{f},
\end{align*}
which is already the left-hand side of \eqref{eq:ofc} in terms of $\tilde{f}$. It only remains to observe that $\mathrm{d}^{\dot{\vec{x}}}_t \tilde{f} = \partial_t f$ and to replace the surface gradient $\nabla_{\cali{M}} \tilde{f}$ by its coordinate expression $J g^{-1} \nabla^2 f$, where $g = J^\top J$ is the coefficient matrix of the Riemannian metric; see e.g.\ \cite{Lee97}.

We highlight the qualitative difference between the constraints \eqref{eq:ofc} and \eqref{eq:ofc2}. Note that in the former the unknown is $\dot{\beta}$, while in the latter it is $\dot{\gamma}_\mathrm{tan} = \dot{\vec{x}}_\mathrm{tan} + J\dot{\beta}$. This means that \eqref{eq:ofc} constrains the motion relative to the tangential surface velocity $\dot{\vec{x}}_\mathrm{tan}$, while \eqref{eq:ofc2} constrains the absolute tangential motion.

The nature of our microscopy data suggests a simple global parametrisation (see Sec.~\ref{sec:num}). We therefore pull the data back to the Euclidean plane and solve \eqref{eq:ofc}. However, equation \eqref{eq:ofc2} is independent of any parametrisation. It can thus serve as a starting point for alternative numerical approaches.

\subsection{Regularisation}\label{sec:regularisation}
From now on we fix an arbitrary $t_0 \in I$ and turn to the actual solution of the parametrised OFC for $(u^1(x), u^2(x))^\top = u(x) = \dot{\beta}(x,t_0)$. Recall that with this notation $u$ contains the coefficients of the tangential vector field $\vec{u} = J\dot{\beta}$ with respect to the tangential basis $(\partial_1\vec{x}, \partial_2\vec{x})$ of $\cali{M}_{t_0}$. Note also that, by fixing $t_0$, there is no more time-dependence in our problem which makes it effectively an optical flow problem on a static surface. Hence we omit any reference to $t_0$ from now on and write $\cali{M}$ instead of $\cali{M}_{t_0}$.

The sought vector field is underdetermined by the OFC alone. We overcome this by minimising a functional that penalises violation of the OFC while imposing an additional smoothness restriction on $\vec{u}$. More precisely, we adopt a recent extension of the original quadratic Horn-Schunck regularisation to a Riemannian setting \cite{LefBai08}. Basically, they propose to minimise
\begin{equation}\label{eq:functional1}
	\cali{E}(u) = \frac{\alpha}{2} \big\| \nabla^2 f \cdot u + \partial_t f \big\|^2_{L^2(\cali{M})} + \frac{1}{2} \big\| Du \big\|^2_{L^2(\cali{M})}.
\end{equation}
Here, $\alpha > 0$ is the regularisation parameter and $Du = (D_ju^i)$ is the $2 \times 2$ matrix containing the coefficient functions of the covariant derivatives
\begin{equation*}
	\nabla_j \vec{u} = \sum_{i=1}^2 D_j u^i \partial_i\vec{x}, \enspace j=1,2,
\end{equation*}
of $\vec{u}$. Using the Christoffel symbols $\Gamma^{i}_{jk}$ (see Sec.~\ref{sec:num}) associated to the parametrisation $\vec{x}$ the coefficients are given by
\begin{equation*}
	D_j u^i  = \partial_j u^i + \sum_{k=1}^2{ \Gamma^{i}_{jk}u^k} , \enspace i,j=1,2.
\end{equation*}
Rewriting \eqref{eq:functional1} as an integral over the coordinate domain, we arrive at the functional
\begin{equation} \label{eq:finalfunctional}
	\cali{E}(u) = \frac{1}{2} \int_\Omega \Big[ \alpha \left( \nabla^2 f  \cdot u + \partial_t f \right)^2 + \norm{Du}^2_F \Big] \sqrt{\det g} \, \mathrm{d}x,
\end{equation}
where $\norm{\cdot}_F$ is the Frobenius norm.
\section{Numerical Solution} \label{sec:num}
We solve the problem of minimising functional $\cali{E}$ via its associated Euler-Lagrange equations. Regarding the integrand of $\cali{E}$ as a function $G(x,u,\nabla^2 u^1, \nabla^2 u^2)$, they read
\begin{align*}
	G_{u^1}	&= \partial_1 G_{\partial_1 u^1} + \partial_2 G_{\partial_2 u^1} \\
	G_{u^2}	&= \partial_1 G_{\partial_1 u^2} + \partial_2 G_{\partial_2 u^2},
\end{align*}
where subscripts of $G$ denote partial derivatives. The resulting pair of linear PDEs is of the form
\begin{equation}\label{eq:pde}
	\begin{aligned}
		\Delta u^1 &= \nabla^2 u^1 \cdot c + \nabla^2 u^2 \cdot d + u \cdot b_1 + a_1 \\
		\Delta u^2 &= \nabla^2 u^2 \cdot c + \nabla^2 u^1 \cdot d + u \cdot b_2 + a_2.
	\end{aligned}
\end{equation}
The coefficient vectors $a,b_1,b_2,c,d$ are rather lengthy functions of the data $f$ and metric tensor $g$, which is why we do not write them out in full here. Letting $\Omega = (0,1)^2$ for simplicity, the natural boundary conditions of the variational problem are
\begin{equation}\label{eq:boundary}
	\partial_j u^i + \sum_k \Gamma^{i}_{jk}u^k = 0, \enspace \mbox{for } x_j \in \{0,1\},
\end{equation}
where $i,j \in \{1,2\}$. In case of a flat manifold, e.g. $\cali{M} = \Omega$, the Euler-Lagrange equations \eqref{eq:pde} reduce to those of the original Horn-Schunck functional and the boundary conditions become the usual homogeneous Neumann ones. For more details on the calculus of variations we refer to \cite{CouHil53}.

Due to the nature of the microscopy data (see Sec.\ \ref{sec:data} and Fig.~\ref{fig:embryo}), the manifold $\cali{M}_t$ modelling the deforming yolk is a surface with boundary that is most easily parametrised as the graph of a function $z: \Omega\times I \to \Reals$. Hence, we set $\vec{x}(x_1,x_2,t) = (x_1,x_2, z(x_1,x_2,t))^\top$. Accordingly, for the metric we get
\begin{equation*}
	g = I_2 + \nabla^2 z \nabla^2 z^\top, \qquad \det g = 1 + \abs{\nabla^2 z}^2,
\end{equation*}
where $I_2\in\Reals^{2 \times 2}$ is the identity matrix. The Christoffel symbols turn out to be
\begin{equation*}
	\Gamma_{jk}^{i} =
	\frac{1}{2} \sum_{m=1}^2{g^{mi} \left( \partial_j g_{km} + \partial_k g_{mj} - \partial_m g_{jk} \right)} =
	\frac{\partial_i z \, \partial_{jk} z}{\det g}.
\end{equation*}
Partial derivatives of $z$ and of the projected data $f$ were approximated by central differences. The system \eqref{eq:pde} with boundary conditions \eqref{eq:boundary} was then solved with a standard finite difference scheme. In the following section numerical results are presented.
\section{Experiments} \label{sec:exp}

\subsection{Data}\label{sec:data}
As mentioned before, the biological motivation for this work are cellular image data of a zebrafish embryo. Endoderm cells expressing green fluorescent protein were recorded via confocal laser-scanning microscopy resulting in time-lapse volumetric (4D) images; see~\cite{MegFra03} for the imaging techniques. This type of image shows a high contrast at cell boundaries and a low signal-to-noise ratio in general. Our videos were obtained during the gastrula period, which is an early stage in the animal's developmental process and takes place approximately five to ten hours post fertilisation. In short, the fish forms on the surface of a spherical-shaped yolk; see e.g.~\cite{KimBalKimUllSchi95} for many illustrations and detailed explanations. For the biological methods such as the fluorescence marker and the embryos used in this work we refer to~\cite{MizVerHeaKurKik08}. The important aspect about endodermal cells is that they are known to form a monolayer during gastrulation~\cite{WarNus99}, meaning that the radial extent is only a single cell. This crucial fact allows for the straightforward extraction of a surface together with a two-dimensional image of the stained cells. Since only a cuboid region of approximately $860 \times 860 \times 340\,\mu\text{m}^{3}$ of the pole region is captured by the microscope, this surface can easily be parametrised; cf. Sec.~\ref{sec:num}. The spatial resolution of the Gaussian filtered images is $512 \times 512$\,pixels and all intensities are given in the interval $[0, 1]$. Our sequence contains $77$ frames recorded in intervals of $240\,\text{s}$ with clearly visible cellular movements and cell divisions.

\subsection{Numerical results} \label{sec:numresults}
In the following we present qualitative results and demonstrate the feasibility of our approach. For every subsequent pair of frames we minimised the functional~\eqref{eq:finalfunctional} as outlined in Sec.~\ref{sec:num}. We chose grid size as well as temporal displacement as $h = 1$ and the regularisation parameter was set to $\alpha = 10$. For demonstration purpose we make use of the standard flow colour-coding~\cite{BakSchaLewRotBla11}, which maps (normalised) flow vectors to a colour space defined inside the unit circle. It is easy to see that the same colours are valid all over the manifold due to the parametrisation.

\begin{figure}[t]
	\centering
	\includegraphics[width=0.49\textwidth]{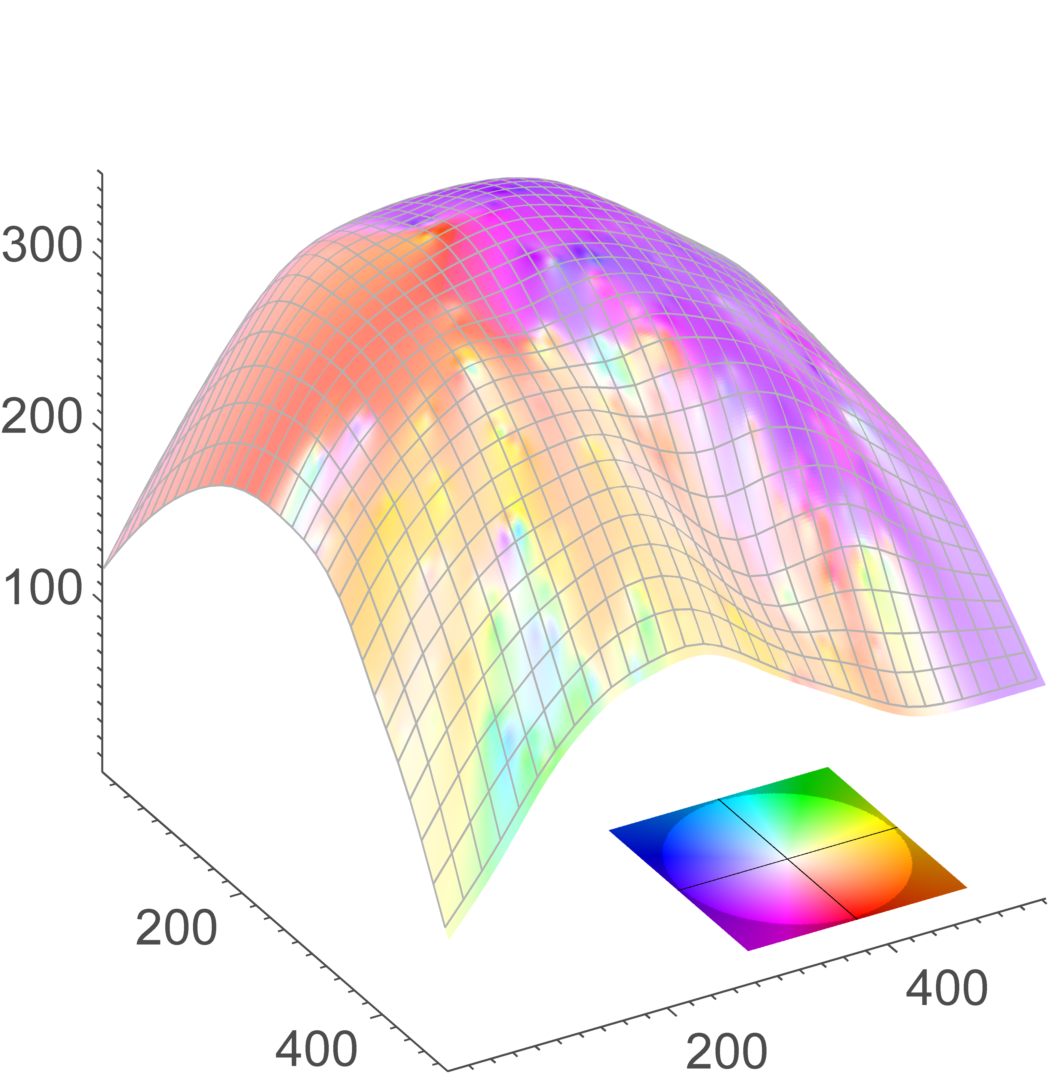}
	\hfill
	\begin{tikzpicture}
		\node[anchor=south west,inner sep=0] at (0,0) {
		\includegraphics[width=0.49\textwidth]{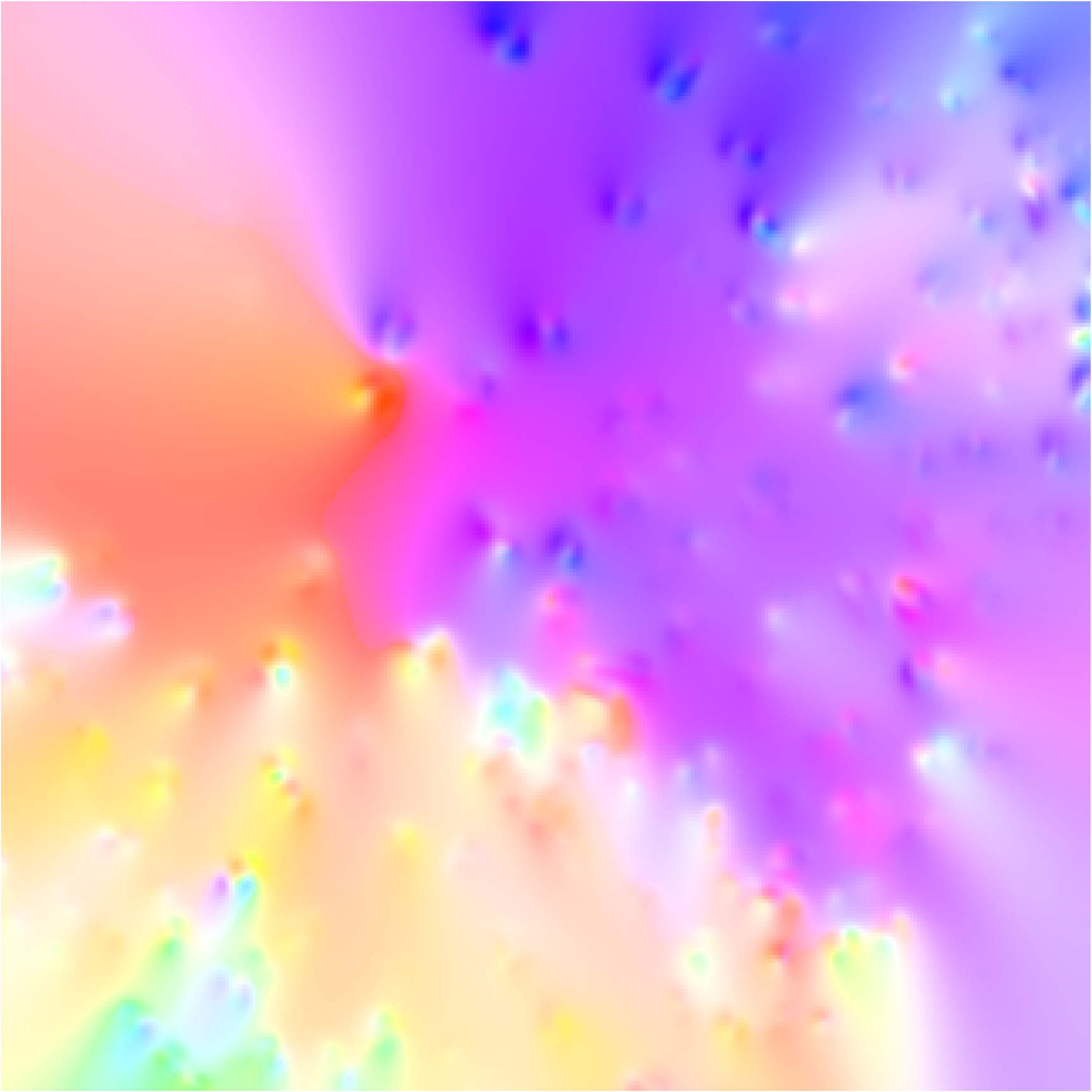}};
	    \draw[red,thick] (0.8cm,0.7cm) rectangle (1.7cm,1.6cm);
	\end{tikzpicture}
	\caption{Optical flow field between frames 57 and 58 of the sequence. Colours indicate direction whereas darkness of a colour indicates the length of the vector. Note that the colour circle has been enlarged for better visibility.}
	\label{fig:flow}
\end{figure}

As representative candidates for this discussion we chose the displacement field between frames 57 and 58 for the following reasons. First, the surface is distinctly developed. Second, a considerable number of cells is present in the image, and third, the interval contains cell divisions. Figure~\ref{fig:flow}, left, shows the colour-coded tangential vector field and the colour space whereas Fig.~\ref{fig:flow}, right, displays the same motion field as computed in the parameter space.\footnote{Some figures may appear in colour only in the online version.} A visual inspection of the dataset shows that cells tend to move towards the embryo's body axis, which roughly runs along the main diagonal in Fig.~\ref{fig:flow}, right. Clearly, the velocity field is sufficiently smooth and suggests this behaviour in an adequate manner on a large scale. The expected change in orientation along the body axis is well represented by the colour shift from orange-yellow below the main diagonal to purplish blue in the region above. On the contrary, the choice of the regularisation parameter ensures that individual movements are well preserved as can be observed from the image. 

\begin{figure}[t]
	\centering
	\includegraphics[width=0.49\textwidth]{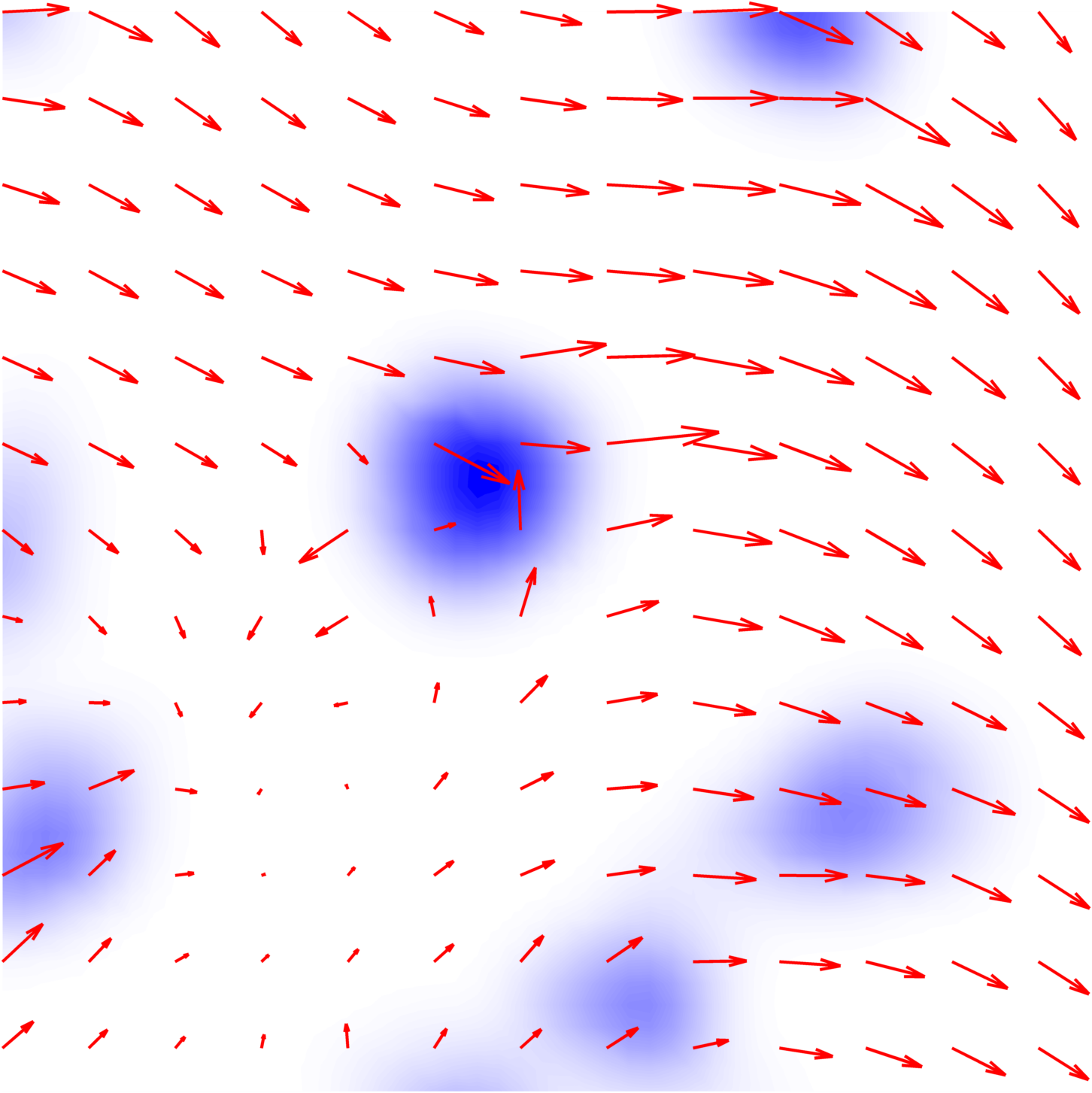}
	\hfill
	\includegraphics[width=0.49\textwidth]{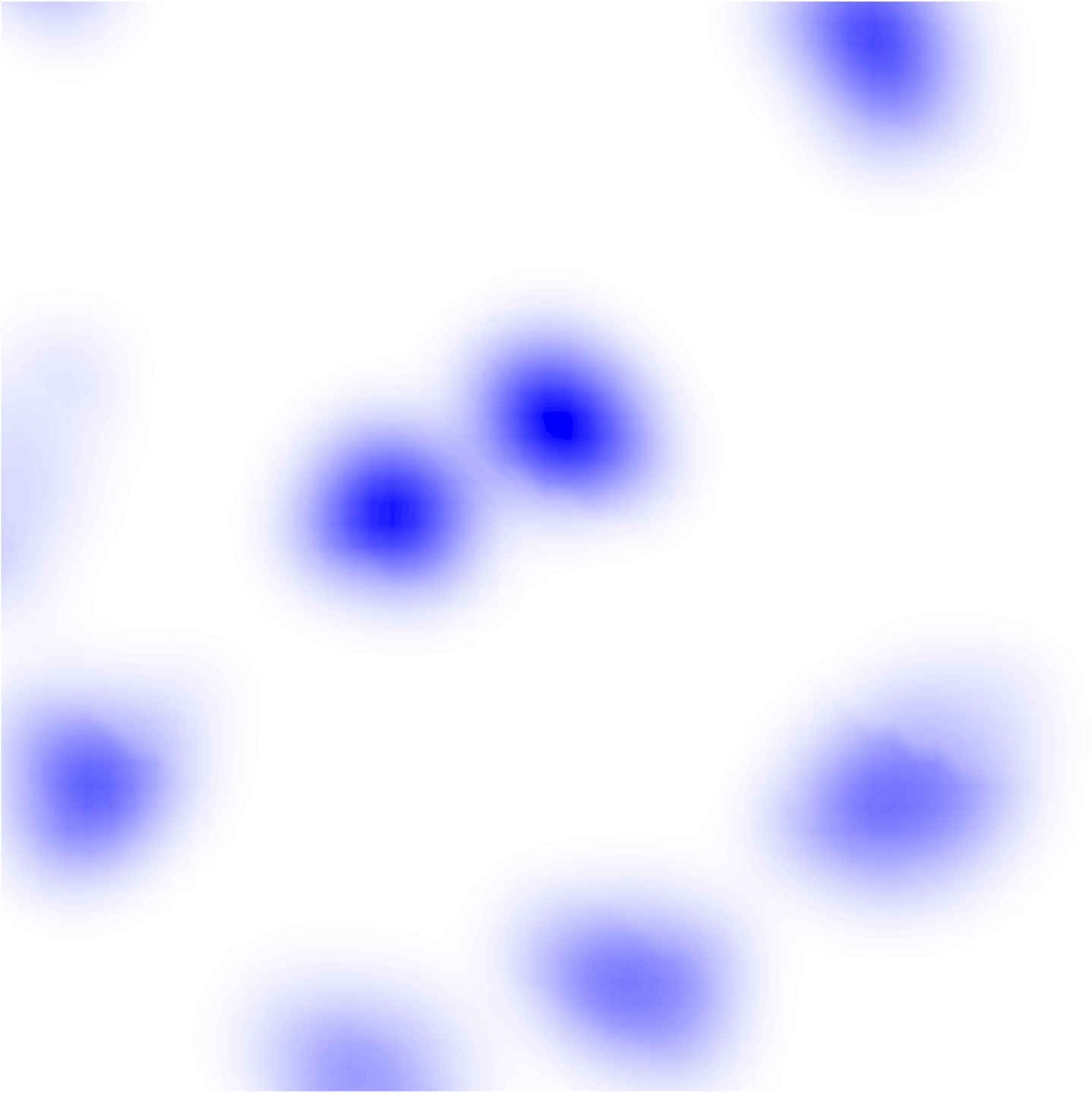}
	\caption{Detailed view of a cell division occurring between frames 57 (left) and 58 (right). All vectors are scaled and only every fourth vector is shown. Intensities are interpolated for smooth illustration}
	\label{fig:division}
\end{figure}

Figure~\ref{fig:division} gives a detailed view of the section outlined by a (red) rectangle in Fig.~\ref{fig:flow}, right. This section was chosen because it depicts a cell division. Figure~\ref{fig:division}, left, and Fig.~\ref{fig:division}, right, display the frames before and after the event, respectively. Moreover, in Fig.~\ref{fig:division}, left, the velocity field is shown. From the raw data we observed that when a cell actually splits, the two daughter cells drift apart in a $180\,^{\circ}$ angle with respect to the mother cell. The displacement field clearly shows the anticipated pattern caused by the diverging daughter cells. In Fig.~\ref{fig:flow}, right, the event is point up by two areas which are coloured mutually opposite with respect to the colour space. Our results suggest that cell division can be indicated reasonably well by our model. Both implementation and data are available on our website.\footnote{\url{http://www.csc.univie.ac.at}}

\section{Conclusion}
Aiming at efficient motion analysis of 4D cellular microscopy data, we generalised the Horn-Schunck method to videos defined on evolving surfaces. The biological fact that the observed cells move along an itself deforming surface allows for motion estimation in 2D (plus time). In the course of this work, we presented two ways to linearise the brightness constancy assumption and showed that one could be obtained from the other and vice versa. The resulting optical flow constraint was solved by means of quadratic regularisation and verified on the basis of the afore-mentioned data. Our qualitative results suggest that both global trends as well as individual movements including cell division are well shown in the surface velocity field. However, so far we only laid the basic groundwork in terms of a mathematical model.

\paragraph{Acknowledgements.}
We thank Pia Aanstad from the University of Innsbruck for sharing her biological insight and for kindly providing the microscopy data. This work has been supported by the Vienna Graduate School in Computational Science (IK I059-N) funded by the University of Vienna. In addition, we acknowledge the support by the Austrian Science Fund (FWF) within the national research networks ``Photoacoustic Imaging in Biology and Medicine" (project S10505-N20, Reconstruction Algorithms for PAI) and ``Geometry + Simulation" (project S11704, Variational Methods for Imaging on Manifolds).

\bibliographystyle{plain}

\def\cprime{$'$} \providecommand{\noopsort}[1]{}\def\cprime{$'$}

\end{document}